\documentclass[sn-mathphys]{sn-jnl}% Math and Physical Sciences Reference Style

\jyear{2021}%

\usepackage{graphicx}
\usepackage{subcaption}
\usepackage{relsize}
\usepackage{enumerate}
\usepackage{soul}

\theoremstyle{thmstyleone}%
\newtheorem{theorem}{Theorem}%  meant for continuous numbers
\newtheorem{lemma}{Lemma}%
\usepackage{thmtools,thm-restate}

\newtheorem{corollary}{Corollary}%

\theoremstyle{thmstyletwo}%
\newtheorem{example}{Example}%

\theoremstyle{thmstylethree}%
\newtheorem{definition}{Definition}%

\begin{document}

\title[The $2$-domination number of cylindrical graphs]{The $2$-domination number of cylindrical graphs}


\author[1,3]{\fnm{Jos\'e Antonio} \sur{Mart\'inez}}\email{jmartine@ual.es}
\equalcont{These authors contributed equally to this work.}

\author[2]{\fnm{Ana Bel\'en} \sur{Casta\~no-Fern\'andez}}\email{acf583@ual.es}
\equalcont{These authors contributed equally to this work.}

\author*[2,3]{\fnm{Mar\'ia Luz} \sur{Puertas}}\email{mpuertas@ual.es}
\equalcont{These authors contributed equally to this work.}

\affil[1]{\orgdiv{Department of Computer Science}, \orgname{Universidad de Almería}, \orgaddress{\street{Carretera Sacramento s/n}, \city{Almer\'ia}, \postcode{04120}, \country{Spain}}}

\affil[2]{\orgdiv{Department of Mathematics}, \orgname{Universidad de Almería}, \orgaddress{\street{Carretera Sacramento s/n}, \city{Almer\'ia}, \postcode{04120}, \country{Spain}}}

\affil[3]{\orgdiv{Agrifood Campus of International Excellence (ceiA3)}, \orgname{Universidad de Almería}, \orgaddress{\street{Carretera Sacramento s/n}, \city{Almer\'ia}, \postcode{04120}, \country{Spain}}}


\abstract{A vertex subset $S$ of a graph $G$ is said to $2$-dominate the graph if each vertex not in $S$ has at least two neighbors in it. As usual, the associated parameter is the minimum cardinal of a $2$-dominating set, which is called the $2$-domination number of the graph $G$. We present both lower and upper bounds of the $2$-domination number of cylinders, which are the Cartesian products of a path and a cycle. These  bounds allow us to compute the exact value of the $2$-domination number of cylinders where the path is arbitrary, and the order of the cycle is $n\equiv 0\pmod 3$ and as large as desired. In the case of the lower bound, we adapt the technique of the wasted domination to this parameter and we use the so-called tropical matrix product to obtain the desired bound. Moreover, we provide a regular patterned construction of a minimum $2$-dominating set in the cylinders having the mentioned cycle order.

%The abstract serves both as a general introduction to the topic and as a brief, non-technical summary of the main results and their implications. Authors are advised to check the author instructions for the journal they are submitting to for word limits and if structural elements like subheadings, citations, or equations are permitted.
}




\maketitle

\section{Introduction}

Graphs are useful tools to model and study problems in networks. For instance, the domination-type parameters can be used to approach the problem of the efficient placement of resources in a network. A \emph{dominating set} in a graph $G$ is a vertex subset $S$ such that each vertex not in $S$ has at least one neighbor in it. The \emph{domination number} of $G$, denoted by $\gamma(G)$, is the cardinal of a minimum dominating set. A general reference of this topic can be found in~\cite{Haynes1998}, where several applications to the domination-type parameters to network problems are shown, such as the social network theory (domination), the selection of sets of representatives (total domination), the planning of bus routes (connected domination), the modeling of computer communication networks (distance-k domination) and the placement of radio stations (broadcast domination).

The variations of the classical definition of domination usually involve additional requirements of either the dominating set or the dominated vertices, which provide properties that are useful to solve a specific network problem. In this paper we focus on the so-called $2$-domination~\cite{Fink1985}, which was originally introduced to further study the applications of domination-type parameters in communication networks under link failures. A \emph{$2$-dominating set} is a vertex subset $S\subseteq V(G)$ such that each vertex not in $S$ has at least two neighbors in it. The \emph{$2$-domination number} $\gamma_2(G)$ is the minimum cardinal of a $2$-dominating set of $G$. The application to the problem of the fault-tolerant placement of sensors in a network was later studied in~\cite{Butjas2018}. Moreover, in the same paper the authors also describe a potential application of the $2$-domination to the data collection problem by a sensor network.

From the computational point of view, it is well known that the decision problem ``Is there a dominating set of the graph $G$ with at most $k$ vertices?'' is NP-complete~\cite{Garey1979}, even if the graph is bipartite or chordal. Similarly, the same problem for $2$-dominating sets is also NP-complete~\cite{Jacobson1989,Bean1994}. This computational complexity makes it interesting to study the $2$-domination number in graph families where the computation of this parameter can be carried out. The Cartesian product graphs are traditionally studied on the domination-type problems since Vizing's conjecture was formulated~\cite{Vizing1963,Vizing1968}. This conjecture, that is still open, proposes a general lower bound for
the domination number of the Cartesian product of two graphs in terms of the domination {numbers} of the factors. A survey on the state of the conjecture can be found in~\cite{Bresar2012}, while recent new approaches are in~\cite{Bresar2021,Gaar2021}. 

The \emph{Cartesian product} of two graphs $G\Box H$ is the graph with vertex set $V(G)\times V(H)$ and such that two vertices $(g_1,h_1), (g_2,h_2)$ are adjacent
in $G\Box H$ if either $g_1=g_2$ and $h_1, h_2$ are adjacent in $H$, or $g_1,g_2$ are adjacent in $G$ and $h_1=h_2$. A general reference about this topic can be found in~\cite{Imrich2000}.  The domination-type parameters in Cartesian product of paths and/or cycles can be computed by using matrix operations~\cite{Klavzar1996,Spalding1998,Guichar2004,Pavlic2012,Goncalves2011,Crevals2014,Rao2019,Carreno2020,Martinez2021}, specifically the  $(\min,+)$ matrix product, also called the tropical product~\cite{Pin1998}. The {\emph{$(\min,+)$ matrix product} $\boxtimes$} is defined over the semi-ring of tropical numbers
$(\mathbb{R}\cup\{\infty\}, \min, +, \infty, 0)$ as $(A \boxtimes B)_{ij}=\min _k(a_{ik}+b_{kj})$. Moreover, for matrix $A$ and $\alpha\in \mathbb{R}\cup\{\infty\}$,
$(\alpha\boxtimes A)_{ij}=\alpha+a_{ij}$. This matrix algebra is well-known and has applications to finite automata~\cite{Pin1998}, statistics~\cite{Omanovic21}, phylogenetics~\cite{Speyer2009}, integer programming~\cite{Butkovic2019}, and other optimization problems~\cite{Krivulin2015}.

The computation of the $2$-domination number of \emph{cylinders}, that is, the Cartesian product of a path and  { a cycle}, has been addressed in cylinders with arbitrary paths and small cycles~\cite{Garzon2022} and in cylinders with small paths and arbitrary cycles~\cite{Garzon2022+}. We now present the computation of this parameter in cylinders with arbitrary paths and with cycles with order $n\equiv 0 \pmod 3$, as large as desired. Following the trends of similar problems for domination-type parameters in selected Cartesian product graphs, we tackle the problem with the aid of specifically developed algorithms, which we implement and run.

This paper is organized as follows. In Section~\ref{section:lower} we obtain both lower and upper bounds of the $2$-domination number in cylinders that allows us to compute this parameter in the mentioned cases. The proof of the key technical result, which in some of its parts has been carried out with the help of a computer, can be found in Section~\ref{section:proof}. Finally, we present our conclusions and analyze the future perspectives of this problem in Section~\ref{section:conclusions}.

\section{A lower bound of the 2-domination number of cylinders}\label{section:lower}

The $2$-domination number of the cylinder $P_m\Box C_n$ is already known in cases $2\leq m\leq 12$ or $3\leq n\leq 15$ (see~\cite{Garzon2022,Garzon2022+}), so throughout this section every cylinder $P_m\Box C_n$ satisfies $m\geq 13$ and $n\geq 16$. Our first purpose is to obtain a lower bound for $\gamma_2(P_m\Box C_n)$ and to this end, we will adapt the technique of the wasted domination proposed in~\cite{Guichar2004} to our {problem. Such} technique has also been used to obtain lower bounds for domination-type parameters in both grids~\cite{Goncalves2011,Rao2019} and cylinders~\cite{Carreno2020,Martinez2021}. We follow these ideas in this paper. 

\begin{figure}[ht]
\centering
\includegraphics[height=0.18\textheight]{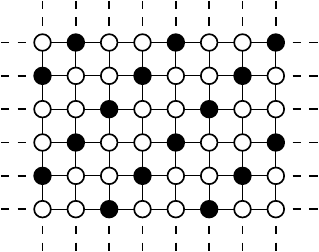}
\caption{The set of black vertices $2$-dominates the infinite grid}
\label{figure:infinite_grid}
\end{figure}

In Figure~\ref{figure:infinite_grid}, we show a $2$-dominating set of an infinite grid consisting of an independent vertex subset $S$ such that every vertex not in it has exactly $2$ neighbors in $S$. This example represents the most efficient way to $2$-dominate the infinite grid because exactly two vertices of $S$ dominate a given vertex not in $S$, in other words, there is no ``wasted'' domination.

\begin{figure}[ht]
\centering
\includegraphics[height=0.35\textheight]{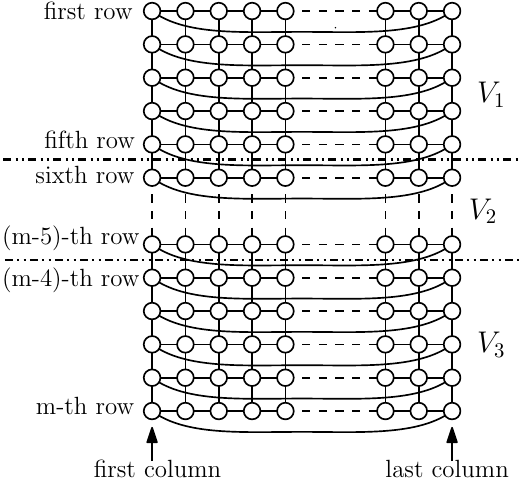}
\caption{Rows and columns of $P_m\Box C_n$}
\label{figure:rows_and_columns}
\end{figure}

However, the cylinder $P_m\Box C_n$ does not admit such regular-patterned $2$-dominating sets and the overlapping of the neighborhoods of the vertices in any $2$-dominating set occurs. Among other reasons, this happens because the borders of the cylinder force the vertices in a $2$-dominating set to stack against each other. Our approach takes into account the behavior of a minimum $2$-dominating set in both borders of the cylinder, and we consider that such borders have five rows (see Figure~\ref{figure:rows_and_columns}). We have carried out some experiments with different sizes of the border and the size of five rows fits our purposes. 
We will use the following notation:
\begin{enumerate}[-]
\item $V=V(P_m\Box C_n)=\{v_{ij}=(a_i,b_j)\colon 1\leq i\leq m, 1\leq j\leq n\}$. 
\par\medskip

\item For $i\in \{1,\dots , m\}$, the {$i$-th} row is the subgraph induced by the vertex subset $\{v_{ij}\colon 1\leq j\leq n\}$ and for $j\in \{1,\dots , n\}$ the {$j$-th} column is the subgraph induced by $\{v_{ij}\colon 1\leq i\leq m\}$.
\par\medskip

\item $V_1=\{v_{ij}\colon 1\leq i\leq 5, 1\leq j\leq n\}$. Moreover, rows from the first to the fourth are called outer rows of $V_1$.
\par\medskip

\item $V_2=\{v_{ij}\colon 6\leq i\leq m-5, 1\leq j\leq n\}$.
\par\medskip

\item $V_3=\{v_{ij}\colon m-4\leq i\leq m, 1\leq j\leq n\}$. Moreover, rows from the {$m$-th} to the {$(m-3)$-th} are the outer rows of $V_3$.
\end{enumerate}
\par\bigskip
It is clear that $V=V_1\cup V_2\cup V_3$ provides a partition of $V$. We call \emph{borders} of $P_m\Box C_n$ to the subgraphs {induced} by $V_1$ and by $V_3$, and the \emph{center} of $P_m\Box C_n$ is the subgraph {induced} by $V_2$. 
\par\bigskip

Let $S$ be a $2$-dominating set of $P_m\Box C_n$. We now consider the partition $S=S_1\cup S_2\cup S_3$, where $S_k=S\cap V_k$ for each $k\in \{1,2,3\}$. For each $k\in \{1,2,3\}$ we define the following vertex subsets:
$$A^S_k=\{v\in V\setminus S \colon v\text{ has at least two neighbors in } S_k\}$$
$$B^S_k=\{v\in V\setminus S \colon v\text{ has exactly one neighbor in } S_k\}$$
\par\medskip

Note that $A^S_k\cap B^S_k=\emptyset$ and $A^S_k\subseteq V_k$ for each $k\in \{1,2,3\}$. In addition, vertices in $B^S_1$ are in the fifth row or in the sixth row of $P_m\Box C_n$, vertices in $B^S_3$ are in the {$(m-4)$-th} row or in the {$(m-5)$-th} row and $B^S_1\cup B^S_3\subseteq A^S_2\cup B^S_2$.

\begin{lemma}\label{lemma:first}
Let $S$ be a $2$-dominating set of $P_m\Box C_n$, then
$$2\left(mn-\vert S\vert \right)\leq  2\left(\vert A^S_1\vert +\vert A^S_2\vert+\vert A^S_3\vert\right) +\vert B^S_1\vert +\vert B^S_2\vert+\vert B^S_3\vert.$$
\end{lemma}

\begin{proof}
On the one hand, $2(\vert A^S_1\vert +\vert A^S_2\vert+\vert A^S_3\vert)$ is twice the number of vertices in $V\setminus S$ having at least two vertices in $S_k$, for some $k$. 

On the other hand, vertices in $V\setminus S$ having exactly one neighbor in $S_1$ and exactly one neighbor in $S_2$ belong to both $B^S_1$ and $B^S_2$ and similarly vertices in $V\setminus S$ having exactly one neighbor in $S_2$ and exactly one neighbor in $S_3$ belong to both $B^S_2$ and $B^S_3$. Therefore, $2(\vert A^S_1\vert +\vert A^S_2\vert+\vert A^S_3\vert) +\vert B^S_1\vert +\vert B^S_2\vert+\vert B^S_3\vert$ is at least twice the number of vertices not belonging to $S$, that is, at least $2(mn-\vert S\vert )$.

Note that neither $A^S_1, B^S_2$ nor $A^S_2, B^S_1$ nor $A^S_3, B^S_2$ nor $A^S_2,B^S_3$ are necessarily disjoint sets, so the inequality is not an equality, in general.
\end{proof}

\begin{lemma}\label{lemma:second}
Let $S$ be a $2$-dominating set of $P_m\Box C_n$ then,
$$4\vert S\vert -2\left(mn-\vert S\vert\right)\geq 4\vert S_1\vert -\left(2\vert A^S_1\vert +\vert B^S_1\vert\right) +4\vert S_3\vert -\left(2\vert A^S_3\vert +\vert B^S_3\vert\right) .$$
\end{lemma}

\begin{proof} 
The partition $S=S_1\cup S_2\cup S_3$ and the inequality obtained in Lemma~\ref{lemma:first} give
$$4\vert S\vert -2(mn-\vert S\vert)=4(\vert S_1\vert +\vert S_2\vert+\vert S_3\vert )-2(mn-\vert S\vert)\geq \sum_{k=1}^{3} {\left(4\vert S_k\vert -(2\vert A^S_k\vert +\vert B^S_k\vert)\right)}.$$ 
So, we just need to prove that $4\vert S_2\vert -(2\vert A^S_2\vert +\vert B^S_2\vert)\geq 0$. Note that $4\vert S_2\vert$ counts the four neighbors of every vertex in $S_2$, including repeats if any. Thus, vertices in $B^S_2$ are being counted and moreover, if a vertex belongs to $A^S_2$, that is,  it is in $V\setminus S$ and it has at least two neighbors in $S_2$ then, it appears twice in such computation. Therefore  $4\vert S_2\vert \geq 2\vert A^S_2\vert +\vert B^S_2\vert$, as desired.
\end{proof}

In order to manage with both borders of the cylinder, we need the following definition. 

\begin{definition}\label{definition:loss}
A vertex subset $R\subseteq V$ is called border-$2$-dominating if it satisfies the following conditions:
\begin{enumerate}[i)]
\item $R\subseteq V_1$ or $R\subseteq V_3$
\item if $R\subseteq V_1$ (respectively $V_3$), {then every} vertex in the outer rows of $V_1$ (respectively $V_3$) has at least two neighbors in $R$,
\item if $R\subseteq V_1$ (respectively $V_3$), {then every} vertex in the {$5$-th}  (respectively the {$(m-4)$-th}) row of $P_m\Box C_n$ has at least one neighbor in $R$.
\end{enumerate}
\end{definition}

Note that every border-$2$-dominating set in $V_3$ is symmetric to another one in $V_1$, so we just focus on $V_1$ and the results are also valid for $V_3$. Let $R$ be a border-$2$-dominating set in $V_1$. We will use the following notation:
$$A^R=\{v\in V\setminus R \colon v\text{ has at least two neighbors in } R\}$$
$$B^R=\{v\in V\setminus R \colon v\text{ {has exactly} one neighbor in } R\}$$

Clearly the vertices in the outer rows of $V_1$ belong to $A^R$ and every vertex in $B^R$, if any, is in the fifth row or in the sixth row of $P_m\Box C_n$. Moreover, if $S\subseteq V$ is a $2$-dominating set{, then} $S_1=S\cap V_1$ (and also $S_3=S\cap V_3$)  is a border-$2$-dominating set.

\begin{definition}
The wasted $2$-domination of a border-$2$-dominating set $R$ of the cylinder $P_m\Box C_n$ is
$$\omega_2(R)=4\vert R\vert -(2\vert A^R\vert +\vert B^R\vert)$$
\end{definition}

For $n\geq 16$ is fixed, the border-$2$-dominating sets of $P_m\Box C_n$ do not depend on $m$ because $V_1$ has five rows whatever $m\geq 13$ is, so we denote 
$$\omega_2(n)=\min \{\omega_2(R)\colon R \text{ is a border-$2$-dominating set of } P_m\Box C_n\}.$$

\begin{restatable}{lemma}{thirdlemma}\label{lemma:third}
Let $n\geq 16$ be an integer. Then
$$\omega_2(n)= \displaystyle\left\{
\begin{array}{ll}
2n+1& \text{if } {n\in \{ 16,19\}}\\
2n& \text{otherwise}
\end{array}
\right.
$$
\end{restatable}

The proof of this key lemma is long and technical and we have carried it out with the help of a computer.  The entire proof can be found in Section~\ref{section:proof}. 

Lemmas~\ref{lemma:second} and~\ref{lemma:third} allow us to obtain the desired lower bound for the $2$-domination number.

\begin{theorem}\label{theorem:bound}
For $m\geq 13$ and $n\geq 16$:
$$\gamma_2(P_m\Box C_n)\geq
\left\{\def\arraystretch{1.5}
\begin{array}{ll}
\displaystyle\frac{(m+2)n+1}{3} & \text{if }  {n\in \{ 16,19\}}\\
\\
\displaystyle\frac{(m+2)n}{3}& \text{otherwise}
\end{array}
\right.
$$
\end{theorem}

\begin{proof}
Let $S$ be a minimum $2$-dominating set of $P_m\Box C_n$. Then, 
$$4\vert S\vert -2(mn-\vert S\vert)=6\vert S\vert -2mn=6 \gamma_2(P_m\Box C_n)-2mn.$$ 
Therefore, by using Lemma~\ref{lemma:second}, we obtain that 
$$6 \gamma_2(P_m\Box C_n)-2mn\geq \omega_2(S_1)+\omega_2(S_3)\geq 2\omega_2(n).$$ 
Finally, by using Lemma~\ref{lemma:third} we obtain the desired inequality 
$$\gamma_2(P_m\Box C_n)\geq \frac{2mn+2\omega_2(n)}{6}=\frac{mn+\omega_2(n)}{3}= \left\{\def\arraystretch{1.5}
\begin{array}{ll}
\displaystyle\frac{(m+2)n+1}{3} & \text{if } {n\in \{ 16,19\}}\\
\\
\displaystyle\frac{(m+2)n}{3}& \text{otherwise}
\end{array}
\right.$$
\end{proof}

The lower bound given in Theorem~\ref{theorem:bound} provides the value of the $2$-domination number if $n\equiv 0 \pmod 3$, as we show in the following corollary.

\begin{corollary}\label{cor:exact}
For $m\geq 8$ and $n\geq 3, n\equiv 0\pmod 3$:
$$\gamma_2(P_m\Box C_n)=\frac{(m+2)n}{3}$$
\end{corollary}

\begin{proof}
Cases with $m=8, 9 ,10, 11, 12$ and any $n\equiv 0\pmod 3$ can be found in \cite{Garzon2022+} and cases with $n=3, 6, 9, 12, 15$ and $m\geq 8$ can be found in \cite{Garzon2022}. 

Assume that $m\geq 13,$ and $n\geq 18, n\equiv 0 \pmod 3$. An upper bound of $\gamma_2(P_m\Box C_n)$ can be obtained by building a specific $2$-dominating set. In Figure~\ref{figure:13_18} we show such set for the smallest cylinder considered now, that is, $P_{13}\Box C_{18}$.

\begin{figure}[ht]
\centering
\includegraphics[height=0.32\textheight]{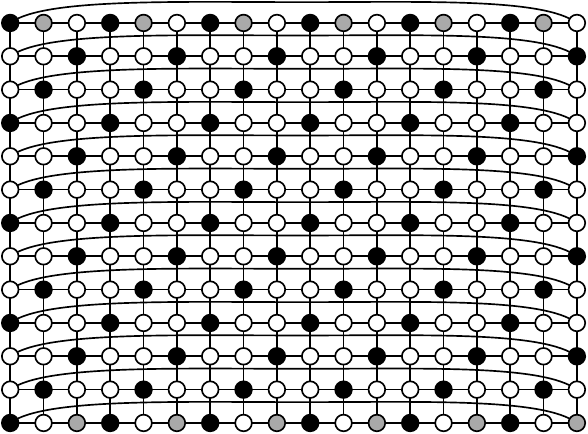}
\caption{Black and gray vertices $2$-dominates $P_{13}\Box C_{18}$}
\label{figure:13_18}
\end{figure}

The general construction for any $m\geq 13$ and $n\geq 18, n\equiv 0 \pmod 3$ is as follows. We pick $\frac{n}{3}$ vertices in each row, the black ones, following the pattern shown in Figure~\ref{figure:13_18}. They $2$-dominate the cylinder except $\frac{n}{3}$ vertices in both the first and the last rows, the gray ones. The resulting set of black and gray vertices $2$-dominates $P_{m}\Box C_{m}$ and it has cardinal
$$ m\frac{n}{3}+ 2\frac{n}{3} =\frac{(m+2)n}{3}\cdot$$ 
Therefore, 
$$\gamma_2(P_m\Box C_n)\leq \frac{(m+2)n}{3}\cdot$$ 
Finally, the lower bound obtained in Theorem~\ref{theorem:bound} gives the desired equality.
\end{proof}

\section{Proof of Lemma~\ref{lemma:third}}\label{section:proof}
We devote this section to the proof of the key lemma that we have used to obtain the lower bound of the $2$-domination number shown in Theorem~\ref{theorem:bound}.
\thirdlemma*

Our approach follows the ideas in \cite{Carreno2020,Guichar2004,Rao2019} and we provide an extension of the technique of the wasted domination, which was first proposed in~\cite{Guichar2004} to compute a lower bound of the domination number in grids. Due to the fact our extension will be applied to cylinders instead of grids, it requires a new tool to obtain the desired lower bound. {Such a tool} is the following theorem from~\cite{Carre1979}, that we quote from~\cite{Klavzar1996} in the version related to the $(\min,+)$ matrix product.

Let $\mathcal{D}$ be a directed graph with vertex set $V(\mathcal{D})= \{w_1, w_2, \dots, w_s\}$ together with a labeling function $\ell$ which assigns an element of the semi-ring $\mathcal{P}=(\mathbb{R}\cup\{\infty\}, \min, +, \infty, 0)$ to every arc of $\mathcal{D}$.
A \emph{path of length $n$} in $\mathcal{D}$ is a sequence of $n$ consecutive arcs $Q=(w_{i_0}w_{i_1})(w_{i_1}w_{i_2})\dots (w_{i_{n-1}}w_{i_n})$ and $Q$ is \emph{a closed path} if $w_{i_0}=w_{i_{n}}$.
Note that the labeling $\ell$ can easily be extended to paths in the following way: $\ell(Q) = \ell(w_{i_0}w_{i_1})+\ell(w_{i_1}w_{i_2})+\dots +\ell(w_{i_{n-1}}w_{i_{n}}).$
\begin{theorem}[\cite{Carre1979,Klavzar1996}]\label{thm:carre}
Let $S_{ij}^n$ be the set of all paths of length $n$ from $w_i$ to $w_j$ in $\mathcal{D}$ and let $A(\mathcal{D})$ be the
matrix defined by
$$A(\mathcal{D})_{ij} =
\left\{
\begin{array}{ll}
\ell(w_i,w_j) & \text{if }(w_i,w_j) \text{ is an arc of } \mathcal{D},\\
\infty & \text{otherwise.}
\end{array}
\right.
$$
If $A(\mathcal{D})^n$ is the {$n$-th} $(\min, +)$ power of $A(\mathcal{D})$, then $(A(\mathcal{D})^n)_{ij}=\min \{\ell(Q)\colon Q\in S_{ij}^n\}.$
\end{theorem}

In order to clarify the steps of the proof of the lemma, a brief sketch is as follows. We define an appropriate digraph $\mathcal{D}$ with an arc labeling $\ell$ such that there is a bijective correspondence between border-$2$-dominating sets $R$ in $P_m\Box C_n$ and closed paths $Q_R$ in $\mathcal{D}$ with length $n$. In addition, this correspondence will provide the equality $\omega_2(R)=\ell(Q_{R})$. Thus, Theorem~\ref{thm:carre} gives: 
\begin{equation*}
\begin{split}
    \omega_2(n) & =\min \{\omega_2(R)\colon R \text{ is a border-$2$-dominating set of } P_m\Box C_n\} \\ 
         & =\min _i(\min \{\ell(Q)\colon Q\in S_{ii}^n\})\\ 
         & =\min_i(A(\mathcal{D})^n)_{ii}
\end{split}
\end{equation*}
Finally, we will compute the necessary $(\min,+)$ powers of the matrix $A(\mathcal{D})$ to obtain the formula for $\omega_2(n)$ claimed in Lemma~\ref{lemma:third}.  

\subsection{Construction of digraph $\mathcal{D}$}

Our next task is to construct an appropriate digraph $\mathcal{D}$. In order to define its vertex set and its arc set, we will need the following definitions.

\begin{definition}\label{def:vertices}

Let $m\geq 13, n\geq 16$ be integers and let $R\subseteq V_1$ be a border $2$-dominating set of $P_m\Box C_n$. The labeling of the vertices of $V_1$ associated to $R$ is the following.

\begin{enumerate}[i)]
    \item $r(v)=0$ if $v\in R$,
    \item $r(v)=1$ if $v\in V_1 \setminus R$ and {$v$} has at least two neighbors in $R$ in its column or in the previous one,
    \item $r(v)=2$ if $v\in V_1 \setminus R$ and {$v$} {has exactly} one neighbor in $R$ in its column or in the previous one,
    \item $r(v)=3$ if $v\in V_1 \setminus R$ and {$v$} has no neighbors in $R$ in its column or in the previous one. 
\end{enumerate}
\end{definition}

The definition of border $2$-dominating set implies that the vertices labeled as $0$ or $1$ can be in any row of $R$. A vertex labeled as $2$, and such that its neighbor in the following column is not in $R$, has exactly one neighbor in $R$ and therefore, it is in the fifth row of $V_1$. Finally, the neighbor in the following column of any vertex with label $3$ has to be in $R$ and moreover, both vertices are in the fifth row of $V_1$.

Let $S$ be a $2$-dominating set then, $R=S\cap V_1$ is a border $2$-dominating set. In this particular case, every vertex with label $2$ such that its neighbor in the following column is not in $R$ must have 
its second neighbor in $S$ in the sixth row of $P_m\Box C_n$. {Similarly, a vertex $v$ with label 3, which is necessarily in the fifth row, has no neighbors in $R=S\cap V_1$ in its column or in the previous one. This means that $v$ has exactly two neighbors in $S$ and moreover, one of them is in its row (the fifth one) and in the following column, while the other one is in the sixth row and the same column as $v$.}

The {$j$-th} column, $1\leq j\leq n$, of any border $2$-dominating set $R$, with vertices $v_{1j}v_{2j}v_{3j}v_{4j}v_{5j}$ can be represented as the $5$-letter word {$\mathbf{p}=r(v_{1j})r(v_{2j})r(v_{3j})r(v_{4j})r(v_{5j})$} in the alphabet $\{0,1,2,3\}$. {Such a word} satisfies some special conditions since $R$ is a border $2$-dominating set: each vertex with label $1$ has at least one neighbor labeled as $0$ in its column, each vertex with label $2$ has at most one neighbor labeled as $0$ in its column and each vertex with label $3$ is in the {fifth} row and its has no neighbor labeled as $0$ in its column. These conditions provide the following definition, that describes the words that can appear as columns of some border $2$-dominating set.

\begin{definition}\label{def:suitable}
A word $\mathbf{p}=p_1 p_2 p_3 p_4 p_5$ of length $5$ in the alphabet $\{0,1,2,3\}$ is called suitable if it satisfies: 
\begin{enumerate}[i)]
    \item $p_k\neq 3$ if $k\in \{1,2,3,4\}$,
    \item $p_1p_2\notin\{ 12, 11\}$,
    \item $p_4p_5\notin\{ 21, 11, 03\}$,
    \item $p_k p_{k+1} p_{k+2}\notin \{020, 111, 112, 211, 212, 113, 213\}$, if $k\in\{1,2,3\}$.
\end{enumerate}
\end{definition}

In a similar way, we can detail rules that explain where two words can appear together as consecutive columns of a border $2$-dominating set. 

\begin{definition}\label{def:follow}
We say that a suitable word $\mathbf{p}=p_1 p_2 p_3 p_4 p_5$ can follow another suitable word $\mathbf{q}=q_1 q_2 q_3 q_4 q_5$ if the following conditions hold.

Conditions for the first letters $p_1$ and $q_1$.
\begin{enumerate}[i)]
    \item if $q_1=0${, then} $p_1=0$ or $p_1=1$ or $\left\{ p_1=2, p_2\neq 0\right\}$,
    \item if $q_1=1${, then} $p_1=0$ or $\left\{p_1=2, p_2=0\right\}$,
    \item if $q_1=2${, then} $p_1=0$.
\end{enumerate}

Conditions for the intermediate letters $p_k$ and $q_k$, with $2\leq k\leq 4$.
\begin{enumerate}[i)]
    \item if $q_k=0${, then} $p_k=0$ or $p_k=1$ or $\left\{p_k=2, p_{k-1}\neq 0, p_{k+1}\neq 0\right\}$,
    \item if $q_k=1${, then} $p_k=0$ or $\left\{p_k=1, p_{k-1}=0, p_{k+1}=0\right\}$ or $\left\{ p_k=2, p_{k-1}=0 \right\}$ or $\left\{p_k=2, p_{k+1}= 0\right\}$,
    \item if $q_k=2${, then} $p_k=0$.
\end{enumerate}

Conditions for the last letters $p_5$ and $q_5$.
\begin{enumerate}[i)]
    \item if $q_5=0${, then} $p_5=0$ or $p_5=1$ or $\left\{p_5=2, p_4\neq 0\right\}$,
    \item if $q_5=1${, then} $p_5=0$ or $\left\{p_5=2, p_4=0\right\}$ or $p_5=3$,
    \item if $q_5=2${, then} $p_5=0$ or $\left\{p_5=2, p_4=0\right\}$ or $p_5=3$,
    \item if $q_5=3${, then} $p_5=0$.
\end{enumerate}
\end{definition}

We now define digraph $\mathcal{D}$ whose vertex set is the set of all suitable words of length $5$ in the alphabet $\{0,1,2,3\}$ and such that there is an arc from the word $\mathbf{q}$ to the word $\mathbf{p}$ if $\mathbf{p}$ can follow $\mathbf{q}$.

\begin{lemma}\label{lem:bijection}
There exists a bijective correspondence between the border $2$-dominating sets of $P_m\Box C_n$ and closed paths of length $n$ in digraph $\mathcal{D}$.
\end{lemma}

\begin{proof}

Let $R\subseteq V_1$ be a border $2$-dominating set of $P_m\Box C_n$ and consider the vertex labeling given by Definition~\ref{def:vertices}. The {$j$-th} column of $R$ is the suitable word $\mathbf{p^j}=r(v_{1j})r(v_{2j})r(v_{3j})r(v_{4j})r(v_{5j})$ and moreover, word  $\mathbf{p^{j+1}}$ can follow word $\mathbf{p^j}$, because  $R$ is border $2$-dominating. So $R$ can be represented by the closed path of length $n$ $Q_R=(\mathbf{p^1} \mathbf{p^2})(\mathbf{p^2} \mathbf{p^3})\dots (\mathbf{p^n}\mathbf{p^1})$.

Conversely, given $Q_R=(\mathbf{p^1} \mathbf{p^2})(\mathbf{p^2} \mathbf{p^3})\dots (\mathbf{p^n}\mathbf{p^1})$, a closed path of length $n$ in the digraph $\mathcal{D}$, we define a vertex subset of $V_1$ in the following way: we identify the {$i$-th} vertex $(1\leq i\leq 5$) of the {$j$-th} column ($1\leq j\leq n$) $v_{ij}$ with the {$i$-th} entry of the {$j$-th} word $p^j_i$ and we define $R_{Q}= \{ v\in V_1\colon v \text{ is identified with } 0\}$. Definition~\ref{def:suitable} and Definition~\ref{def:follow} gives that each vertex in the outer rows of $V_1$ has at least two neighbors in $R_{Q}$, and each vertex in the {fifth} row of $V_1$ has at least one neighbor in $R_{Q}$. Therefore, $R_{Q}$ is a border $2$-dominating set associated to the closed path $Q$.

This correspondence between border $2$-dominating sets of $P_m\Box C_n$ and closed path of length $n$ of $\mathcal{D}$ is the desired bijection.
\end{proof}

\subsection{Definition of the arc labeling $\ell$}

We now focus on the definition of the arc labeling $\ell$ of $\mathcal{D}$ such that, for each border $2$-dominating set $R$,
$$\ell (Q_R)=\omega_2(R)=4\vert R\vert -(2\vert A^R\vert +\vert B^R\vert)$$ 
where $A^R=\{v\in V\setminus R \colon v\text{ has at least two neighbors in } R\}$
and $B^R=\{v\in V\setminus R \colon $ $ v\text{ has at exactly one neighbor in } R\}$. 

To this end, firstly note that the labels of the vertices of $\mathcal{D}$ given in Definition~\ref{def:vertices}, easily gives that
\begin{equation*}
\begin{split}
    R= &\{v\in V_1\colon r(v)=0 \}\\
    A^R  = &\{ v\in V_1\colon r(v)=1 \}\ \mathlarger{\cup}\ \{ v\in V_1\colon r(v)=2, \text{ and its neighbor in the}\\
           &  \text{following column has label } 0 \}\\
    B^R  = &\{ v\in V_1\colon r(v)=3, \text{ and it is in the {fifth} row} \}  \ \mathlarger{\cup}\  \{ v\in V_1\colon r(v)=2, \\
    &\text{it is in the {fifth} row, and its neighbor in the following column is not }\\
             & \text{labeled as  } 0\}\ \mathlarger{\cup}\  \{ u\in V\colon v \text{ is in the sixth row, and its neighbor in the}\\
           & \text{{fifth} row has label } 0 \} 
\end{split}
\end{equation*}

To compute the cardinal of the sets $R, A^R, B^R$ by using the arcs of the closed path $Q_R$, we follow the ideas and the notation of~\cite{Guichar2004}, by adapting the concept of newly dominated vertices to our case. 

\begin{definition}\label{def:newly}
Let $\mathbf{p}$ and $\mathbf{q}$ two suitable words such that $\mathbf{p}$ can follow $\mathbf{q}$. 
\begin{enumerate}[i)]
    \item The dominating vertices of the arc $(\mathbf{q}\mathbf{p})$ are the vertices in $\mathbf{p}$ with label $0$. We define $d(\mathbf{q} \mathbf{p})$ as the number of dominating vertices of the arc.

\item  The newly $2$-dominated vertices of the arc $(\mathbf{q}\mathbf{p})$ are the vertices in $\mathbf{p}$ with label $1$ and the vertices in $\mathbf{q}$ with label $2$ whose neighbor in $\mathbf{p}$ is labeled with $0$. We define $nd_2(\mathbf{q} \mathbf{p})$ as the number of newly $2$-dominated vertices of the arc.

\item  The newly $1$-dominated vertices of the arc $(\mathbf{q}\mathbf{p})$ are the vertices in the fifth row of $\mathbf{q}$ with label $3$, the vertices of the fifth row of $\mathbf{q}$ with label $2$ whose neighbor in $\mathbf{p}$ is not labeled with $0$ and the neighbor in the sixth row of the cylinder $P_m\Box C_n$ of each vertex in the fifth row of $\mathbf{p}$ with label $0$. We define $nd_1(\mathbf{q} \mathbf{p})$ as the number of newly $1$-dominated vertices of the arc.
\end{enumerate}
\end{definition}

We illustrate these definitions with the following example.

\begin{example}
Consider the suitable words $\mathbf{q}=21013$, $\mathbf{p}=01010$. Note that $\mathbf{p}$ can follow  $\mathbf{q}$, so $(\mathbf{q}\mathbf{p})$ is an arc in $\mathcal{D}$. We draw them as columns in Figure~\ref{fig:example_a}. 

Moreover, in Figure~\ref{fig:example_b} underlined vertices are the dominating vertices of the arc, vertices with a square are the newly $2$-dominated and vertices with a circle are the newly $1$-dominated. 
Note that the newly $1$-dominated vertex $u$ is in the sixth row. 

In addition, in this case $d(\mathbf{q}\mathbf{p})=3, nd_2(\mathbf{q}\mathbf{p})=3$ and $nd_1(\mathbf{q}\mathbf{p})=2$.

\begin{figure}[ht]
     \centering
     \begin{subfigure}{.21\textwidth}
     \centering
        \subcaptionbox{{The arc $(\mathbf{q} \mathbf{p})$ represented as two consecutive columns}\label{fig:example_a}}[3.5cm] {\includegraphics[width=0.5\textwidth]{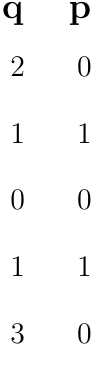}}
     \end{subfigure}
     \hspace{2cm}
     \begin{subfigure}{.2\textwidth}
     \centering
         \subcaptionbox{{Dominating, newly $2$-dominated and newly $1$-dominated vertices} \label{fig:example_b}}[3.5cm]
         {\includegraphics[width=0.5\textwidth]{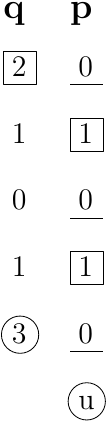}}
        
     \end{subfigure}
        \caption{Newly dominated vertices}
        \label{fig:example}
\end{figure}

\end{example}

We now complete the construction of the digraph $\mathcal{D}$ with the definition of the arc labeling $\ell$.

\begin{definition}\label{def:labeling}
Let $\mathbf{q}\mathbf{p}$ an arc of digraph $\mathcal{D}$. Then 
$$\ell(\mathbf{q}\mathbf{p})=4d(\mathbf{q} \mathbf{p})- \left(2nd_2(\mathbf{q} \mathbf{p})+nd_1(\mathbf{q} \mathbf{p})\right).$$
\end{definition}

In the following lemma we relate the arc labeling with the wasted $2$-domination. 

\begin{lemma}
Let $R$ be a border $2$-dominating set of $P_m\Box C_n$ and let $Q_R$ be its associated closed path in $\mathcal{D}$. Then,
$$\ell (Q_R)=\omega_2(R)=4\vert R\vert -(2\vert A^R\vert +\vert B^R\vert).$$
\end{lemma}

\begin{proof}
Denote $Q_R=(\mathbf{p^1} \mathbf{p^2})(\mathbf{p^2} \mathbf{p^3})\dots (\mathbf{p^n}\mathbf{p^1})$. Firstly, for each arc in $Q_R$, Lemma~\ref{lem:bijection} and Definition~\ref{def:newly} give that the dominating vertices belong to $R$ (they have label $0$), the newly $2$-dominated vertices belong to $A^R$ (they have at least $2$ neighbors in $R$) and the newly $1$-dominated vertices belong to $B^R$ (they have a unique neighbor in $R$). 

Conversely, every vertex in $R$ is a dominating vertex of a unique arc in $Q_R$, every vertex in $A^R$ is a newly $2$-dominated vertex of a unique arc in $Q_R$ and every vertex in $B^R$ is a newly $1$-dominated vertex of a unique arc in $Q_R$. Therefore, by using the definition of the labeling $\ell$, we obtain that

\begin{equation*}
\begin{split}
    \ell(Q_R)  = & \sum_{j=1}^{n-1} \ell ( \mathbf{p^j} \mathbf{p^{j+1}} ) + \ell ( \mathbf{p^n} \mathbf{p^1} ) \\
               = & \sum_{j=1}^{n-1} 4 d( \mathbf{p^j} \mathbf{p^{j+1}} )- \left( 2 nd_2( \mathbf{p^j} \mathbf{p^{j+1}} ) + nd_1 ( \mathbf{p^j} \mathbf{p^{j+1}} ) \right) + \\
                 & 4 d( \mathbf{p^n} \mathbf{p^{1}} )- \left( 2 nd_2( \mathbf{p^n} \mathbf{p^{1}} ) + nd_1 ( \mathbf{p^n} \mathbf{p^{1}} ) \right)\\
               = & \sum_{j=1}^{n-1} 4 d( \mathbf{p^j} \mathbf{p^{j+1}} ) + 4  d( \mathbf{p^n} \mathbf{p^{1}} ) +\\
                & (-2) \left( \sum_{j=1}^{n-1} nd_2( \mathbf{p^j} \mathbf{p^{j+1}}) + nd_2( \mathbf{p^n} \mathbf{p^{1}})\right) +\\
                & (-1) \left( \sum_{j=1}^{n-1} nd_1( \mathbf{p^j} \mathbf{p^{j+1}}) + nd_1( \mathbf{p^n} \mathbf{p^{1}})\right)\\
               =& 4\vert R\vert -(2\vert A^R\vert +\vert B^R\vert)=\omega_2(R).    
\end{split}
\end{equation*}
\end{proof}

\subsection{The $(\min ,+)$ powers of the matrix  $A(\mathcal{D})$  }
We now prove the result that will allow us to compute $\omega_2(n)$ by using the $(\min,+)$ powers of the matrix $A(\mathcal{D})$.

\begin{lemma}\label{lem:powers}
Let $m\geq 13, n\geq 16$ be integers. Let $\mathcal{D}$ be the digraph with the arc labeling $\ell$ that we have constructed above. Define the matrix
$$A(\mathcal{D})_{\mathbf{q}\mathbf{p}} =
\left\{
\begin{array}{ll}
\ell(\mathbf{q}\mathbf{p}) & \text{if }(\mathbf{q}\mathbf{p}) \text{ is an arc of } \mathcal{D},\\
\infty & \text{otherwise.}
\end{array}
\right.
$$
Then, 
\begin{equation*}
\begin{split}
\omega_2(n)=&\min \{ \omega_2(R)\colon R \text{ is a border-2-dominating set of } P_m\Box C_n\}\\
=& \min_{\mathbf{p}\in V(\mathcal{D}) } (A(\mathcal{D})^n)_{\mathbf{p}\mathbf{p}}
\end{split}
\end{equation*}
\end{lemma}

\begin{proof}
Firstly, Theorem~\ref{thm:carre} and Lemma~\ref{lem:powers} give that 
\begin{equation*}
\begin{split}
(A(\mathcal{D})^n)_{\mathbf{p}\mathbf{p}}= & \min \{\ell(Q)\colon Q\in S_{\mathbf{p}\mathbf{p}}^n\}\\
= &  \min \{\ell(Q)\colon Q \text{ closed path with length n, from } \mathbf{p} \text{ to } \mathbf{p}\} \\
= & \min \{ \omega_2(R)\colon R \text{ border 2-dominating set with } \mathbf{p} \text{ as a first column}\}
\end{split}
\end{equation*}
Therefore 
\begin{equation*}
\begin{split}
\min_{\mathbf{p}} (A(\mathcal{D})^n)_{\mathbf{p}\mathbf{p}}= & \min_{\mathbf{p}} \{ \omega_2(R)\colon R \text{ border 2-dominating set with } \mathbf{p} \text{ as a first column}\}\\
 =& 
\min \{ \omega_2(R)\colon R \text{ border 2-dominating set}\}\\
= &\omega_2(n)
\end{split}
\end{equation*}
\end{proof}

We will use the following standard property of the  $(\min,+)$ matrix product, whose proof we include for the shake of completeness.

\begin{lemma}\label{lem:recurrence}
Let $M$ be a square matrix. Suppose that there exist natural numbers $n_0,a,b$ such that $M^{n_0+a}=b\boxtimes M^{n_0}$. Then, $M^{n+a}=b\boxtimes M^{n}$, for every $n\geq n_0$.
\end{lemma}

\begin{proof}
By hypothesis, $M^{n_0+a}=b\boxtimes M^{n_0}$. Let $n\geq n_0$ be such that $M^{n+a}=b\boxtimes M^{n}$ then,
$M^{(n+1)+a}=M\boxtimes M^{n+a}=M\boxtimes(b\boxtimes M^{n})=b\boxtimes(M\boxtimes M^{n})=b\boxtimes M^{n+1}$. 
\end{proof}

The application of this property of the $(\min,+)$ matrix product provides the following result.

\begin{lemma}\label{lem:final}
If there exist natural numbers $n_0,a,b$ such that $A(\mathcal{D})^{n_0+a}=b\boxtimes A(\mathcal{D})^{n_0}${, then} $\omega_2(n+a)-\omega_2(n)=b$, for every $n\geq n_0$.
\end{lemma}

\begin{proof}
By using Lemma~\ref{lem:recurrence}, we obtain that $A(\mathcal{D})^{n+a}=b\boxtimes A(\mathcal{D})^{n}$, for every $n\geq n_0$. Therefore, by Lemma~\ref{lem:powers}, 
\begin{equation*}
\begin{split}
\omega_2(n+a)= & \min_{\mathbf{p}} (A(\mathcal{D})^{n+a})_{\mathbf{p}\mathbf{p}}\\
 = &\min_{\mathbf{p}} (b\boxtimes A(\mathcal{D})^{n})_{\mathbf{p}\mathbf{p}}\\
 =& b+\min_{\mathbf{p}} (A(\mathcal{D})^{n})_{\mathbf{p}\mathbf{p}}\\
 = &b+\omega_2(n)
\end{split}
\end{equation*}
\end{proof}

\subsection{Proof of the lemma}
The finite difference equation given by Lemma~\ref{lem:final} has as a solution the desired formula for $\omega_2(n)$. While it is true that the existence of the natural numbers $n_0,a,b$ of the hypothesis of Lemma~\ref{lem:final} is not theoretically ensured, we have found them by effectively computing the matrix $A(\mathcal{D})$ associated to digraph $\mathcal{D}$ with arc labeling $\ell$, and computing enough $(\min,+)$ powers of such matrix. The size of the matrix makes it necessary to use a computer to perform the computations, and to this end we have designed the following algorithm.

\begin{algorithm}[H]
  \caption{Computation of the matrix $A(\mathcal{D})$ and its powers}
  \begin{algorithmic}[1]
    \Ensure{$A(\mathcal{D})^k$ and  $\min_{\mathbf{p}\in V(\mathcal{D}) } (A(\mathcal{D})^k)_{\mathbf{p}\mathbf{p}}$ for $k$ big enough}
    \State Compute the suitable words\Comment{Definition~\ref{def:suitable}}
    \State Compute the matrix $A(\mathcal{D})$ \Comment{Definitions~\ref{def:newly} and~\ref{def:labeling} and Lemma~\ref{lem:powers}}
    \State Compute $A(\mathcal{D})^k$, for $k$ big enough\Comment{$(\min ,+)$ matrix product}
    \State Compute $\min_{\mathbf{p}} (A(\mathcal{D})^k)_{\mathbf{p}\mathbf{p}}$, for each $k$
  \end{algorithmic}
\label{code:matrix}
\end{algorithm}

\thirdlemma*
\begin{proof}

In Step 1 of Algorithm~\ref{code:matrix}, we first obtain the set of suitable words by computing the $5$-element variations with repetition of the elements in the alphabet $\{0,1,2,3\}$ and keeping those of them that satisfy the conditions given in Definition~\ref{def:suitable}. There are $625(=5^4)$ words of which $111$ are suitable words. Therefore, digraph $\mathcal{D}$ has $111$ vertices.

In Step 2, we compute the labeling of the arcs of  digraph $\mathcal{D}$, by using Definitions~\ref{def:newly} and~\ref{def:labeling}, and the matrix $A(\mathcal{D})$ as described in Lemma~\ref{lem:powers}. The matrix has $111$ rows and columns.

In Step 3, we {compute} successive $(\min,+)$ powers of the matrix $A(\mathcal{D})$. We {compare} them to each other until we {obtain} that $A(\mathcal{D})^{46}=2 \boxtimes A(\mathcal{D})^{45}$. Therefore, $n_0=45,a=1,b=2$ are the values to apply Lemma~\ref{lem:final}, that gives the following finite difference equation for $\omega_2(n)$:
$$\omega_2(n+1)-\omega_2(n)=2, \text{ for every } n\geq 45.$$ 

In Step 4, the minimum of the main diagonal of each power of $A(\mathcal{D})$ {is obtained}. In particular, $\omega_2(45)=\min_{\mathbf{p}} (A(\mathcal{D})^{45})_{\mathbf{p}\mathbf{p}} =90$ and therefore, the solution of the finite difference equation is 
$$\omega_2(n)=2n, \text{ for every } n\geq 45.$$ 

In addition, for $16\leq n\leq 44$, $\omega_2(16)=\min_{\mathbf{p}} (A(\mathcal{D})^{16})_{\mathbf{p}\mathbf{p}}=33, \omega_2(19)=\min_{\mathbf{p}} (A(\mathcal{D})^{19})_{\mathbf{p}\mathbf{p}}=39$ and $\omega_2(n)=\min_{\mathbf{p}} (A(\mathcal{D})^{n})_{\mathbf{p}\mathbf{p}}=2n$, otherwise. This concludes the proof.
\end{proof}

The implementation of the algorithm has been developed in C language and we have run it on an Intel(R) Xeon(R) CPU E5-2650 at 2.00GHz processor, with 64 GB of memory. In spite of the algorithmic complexity of the matrix product, the running time of Algorithm~\ref{code:matrix} is less than $1$ second because the matrix size, $111$ rows and columns, is small for a computational approach.

\section{Conclusions}\label{section:conclusions}
We have computed the $2$-domination number of $P_m\Box C_n$ for $m\geq 13$ and $n\geq 16, n\equiv 0\pmod 3$, obtaining the formula $\gamma_2(P_m\Box C_n)=\frac{(m+2)n}{3}$, which is the same as the small cases obtained in~\cite{Garzon2022,Garzon2022+}.

The technique we have used is a modification of the wasted domination and the newly dominated vertices introduced in~\cite{Guichar2004} to obtain a lower bound of the domination number of grids, that is, the Cartesian product of two paths. Such technique was later applied to compute the exact value of the  domination number of grids and also other domination-type parameters in both grids and cylinders. We have adapted it to the case of $2$-domination in cylinders, and these ideas are a cornerstone of the computation of domination-type parameters in the Cartesian products of paths and/or cycles.

We have used borders of the cylinder with five rows because this size suits our purposes, but we have tested other sizes. For borders with $3$ or $4$ rows, smaller lower bounds can be obtained, so they will not give the desired exact value. Even in these cases, the matrices involved are too large to manage them without the aid of a computer, given that they have $18$ and $45$ rows and columns respectively. We have also tested larger border sizes with $6,7$ and $8$ rows and the matrix size increases to $276$, $687$ and $1707$ rows and columns respectively. Moreover, the algorithm running times also increase to $10$ seconds, $2$ minutes, and $50$ minutes respectively, so they increase quicker than the matrix size. This means that much larger border sizes could not be tested due to the complexity of the problem.

{Obviously, if the border has $6, 7$ or $8$ rows and $n\equiv 0\pmod 3$, then the lower bound agrees with the bound obtained in Theorem}~\ref{theorem:bound}{, because it gives the exact value. However, if $n\not\equiv 0\pmod 3$, then} the bound slightly increases, which leads us to the following thought. In the case of grids, it is expected that domination-type parameters show a non-regular behavior when both paths are small, but a unique formula appears for large enough paths, as happens for the domination number~\cite{Goncalves2011}, the Roman domination number, and the $2$-domination number~\cite{Rao2019}. The study of cylinders is not completed and domination number~\cite{Carreno2020} and the Roman domination number~\cite{Martinez2021} have been computed for small cases and for cylinders with cycles with non-bounded particular orders. It seems that there will be no single formula for the large cases, but different ones depending on the parity of the cycle order. 

The cause of this different behavior of grids and cylinders can be described in terms of the wasted domination, for the corresponding parameter. In big grids, the center is expected to have no wasted domination and the stacking of the vertices of a minimum dominating-type set just occurs next to the borders. However, in cylinders such stacking takes place all around the graph except in the cases that the cycle order allows to uniformly distribute the vertices of such sets. For this reason, we think that computing the wasted domination in the borders of the cylinder will not be enough in all cases, and additional techniques will be necessary to completely compute the domination number and other domination-type parameters in cylinders.

\section*{Funding} This work was partially supported by the grants of the Spanish Ministry of Science and Innovation PID2021-123278OB-I00 and PID2019-104129GB-I00/AEI/10.13039/501100011033, both
funded by MCIN/AEI/10.13039/ 501100011033 and by ``ERDF A way to make Europe.''

\bibliography{bibfile}

\end{document}